\documentclass[a4paper]{amsart} 
\usepackage{amsfonts}
\title[The sigma-sequence and occurrences of some patterns, subsequences and subwords]{The sigma-sequence and counting occurrences of some patterns, subsequences and subwords}
\author{Sergey Kitaev}
\email{kitaev@math.chalmers.se} 
\address{Matematik, Chalmers tekniska h\"ogskola och G\"oteborgs
universitet, S-412~96  G\"oteborg, Sweden}

\newtheorem{prop}{Proposition}
\newtheorem{lemma}[prop]{Lemma}
\newtheorem{cor}[prop]{Corollary}
\newtheorem{thm}[prop]{Theorem}
\theoremstyle{definition}

\newtheorem{ex}[prop]{Example}
\newtheorem{remark}[prop]{Remark}

\def\mn{\mbox{-}}
\def\SS{{\mathcal S}}
\def\KK{{\mathcal K}}

\begin{document}
\begin{abstract} 
We consider {\em sigma-words}, which are words used by Evdokimov in the construction of the sigma-sequence \cite{Evdok}. We then find the number of occurrences of certain patterns and subwords in these words.
\end{abstract}

\maketitle 

\thispagestyle{empty}

\section{Introduction and Background}

We write permutations as words $\pi=a_1 a_2\cdots a_n$, whose letters are distinct and usually consist of the integers $1,2,\ldots,n$.

An occurrence of a pattern $p$ in a permutation $\pi$ is ``classically'' defined as a subsequence in $\pi$ (of the same length as the length of $p$) whose letters are in the same relative order as those in $p$. Formally speaking, for $r \leq n$, we say that a permutation $\sigma$ in the symmetric group ${\mathcal S}_n$ has an occurrence of the pattern $p \in {\mathcal S}_r$ if there exist $1 \leq i_1 < i_2 < \cdots < i_r \leq n$ such that $p = \sigma(i_1)\sigma(i_2) \ldots \sigma(i_r)$ in reduced form. The {\em reduced form} of a permutation $\sigma$ on a set $\{ j_1,j_2, \ldots ,j_r \}$, where $j_1 < j_2 < \cdots <j_r$, is a permutation ${\sigma}_1$ obtained by renaming the letters of the permutation $\sigma$ so that $j_i$ is renamed $i$ for all $i \in \{1, \ldots ,r\}$. For example, the reduced form of the permutation 3651 is 2431. The first case of classical patterns studied was that of permutations avoiding a pattern of length 3 in $\SS_3$. Knuth \cite{Knuth} found that, for any $\tau\in\SS_3$, the number $|\SS_n(\tau)|$ of $n$-permutations avoiding $\tau$ is $C_n$, the $n$th Catalan number.
Later, Simion and Schmidt \cite{SimSch} determined the number $|\SS_n(P)|$ of permutations in $\SS_n$ simultaneously avoiding any given set of patterns $P\subseteq\SS_3$.

In \cite{BabStein} Babson and Steingr\'{\i}msson introduced {\em generalised permutation patterns} that allow the requirement that two adjacent letters in a pattern must be adjacent in the permutation. In order to avoid confusion we write a "classical" pattern, say $231$, as $2\mn 3\mn 1$, and if we write, say $2\mn 31$, then we mean that if this pattern occurs in the permutation, then the letters in the permutation that correspond to $3$ and $1$ are adjacent. For example, the permutation $\pi=516423$ has only one occurrence of the pattern $2\mn 31$, namely the subword 564, whereas the pattern $2\mn 3\mn 1$ occurs, in addition, in the subwords 562 and 563. A motivation for introducing these patterns in \cite{BabStein} was the study of Mahonian statistics. A number of interesting results on generalised patterns were obtained in \cite{Claes}. Relations to several well studied combinatorial structures, such as set partitions, Dyck paths, Motzkin paths and involutions, were shown there.

Burstein \cite{Burstein} considered words instead of permutations. In particular, he found the number $|[k]^n(P)|$ of words of length $n$ in a $k$-letter alphabet that avoid all patterns from a set $P\subseteq\SS_3$ simultaneously. Burstein and Mansour \cite{BurMans1} (resp. \cite{BurMans2,BurMans3}) considered forbidden patterns (resp. generalized patterns) with repeated letters.

The most attention, in the papers on classical or generalized patterns, is paid to finding exact formulas and/or generating functions for the number of words or permutations avoiding, or having $k$ occurrences of, certain patterns. In~\cite{KitMans1} the authors suggested another problem, namely counting the number of occurrences of certain patterns in certain words. These words were chosen to be the set of all finite approximations of a sequence generated by a {\em morphism} with certain restrictions. A motivation for this choice was the interest in studying classes of sequences and words that are defined by iterative schemes~\cite{Lothaire,Salomaa}. In~\cite{KitMans2} the authors also studied the number of occurrences of certain patterns in certain words. But there they choose these words to be the subdivision stages from which the {\em Peano curve} is obtained. The authors called these words the {\em Peano words}. The Peano curve was studied by the Italian mathematician Giuseppe Peano in 1890 as an example of a continuous space filling curve. 

In the present paper we consider the {\em sigma-words}, which are words used by Evdokimov in construction of the {\em sigma-sequence} \cite{Evdok}. Evdokimov used this sequence to construct chains of maximal length in the $n$-dimensional unit cube. Independent interest to the sigma-sequence appears in connection with the well-known {\em Dragon curve}, discovered by physicist John E. Heighway and defined as follows: we fold a sheet of paper in half, then fold in half again, and again, etc. and then unfold in such way that each crease created by the folding process is opened out into a 90-degree angle. The ``curve'' refers to the shape of the partially unfolded paper as seen edge on. If one travels along the curve, some of the creases will represent turns to the left and others turns to the right. Now if 1 indicates a turn to the right, and 2 to the left, and we start travelling along the curve indicating the turns, we get the sigma-sequence~\cite{Evdokimov}. In \cite{Kitaev} the sigma-sequence was studied from another point of view. It was proved there that this sequence cannot be defined by iterated morphism. 

Since the sigma-sequence $w_{\sigma}$ is a sequence in a 2-letter alphabet, we consider patterns in 2-letter alphabets. Moreover, the patterns in a 1-letter alphabet (for example $1\mn 1\mn 1$) correspond to two subsequences (for this example, these subsequences are $1\mn 1\mn 1$ and $2\mn 2\mn 2$), whereas the patterns in a 2-letter alphabet (with at least one letter 2) uniquely determine the subsequences in $w_{\sigma}$ that correspond to them, and conversely. For example, an occurrence of the pattern $1\mn 2\mn 1$ is an occurrence of the subsequence $1\mn 2\mn 1$, whereas an occurrence of the subsequence (subword) $211$ is an occurrence of the pattern $211$. Thus, any our result for a pattern, can be interpreted in term of subsequences or subwords, depending on the context, and conversely.

In our paper we give either an explicit formula or recurrence relation for the number of occurrences for some classes of patterns, subwords and subsequences in the sigma-words. In particular, Theorem~\ref{thm1}, allows to find the number of occurrences of an arbitrary generalized pattern without internal dashes of length $\ell$, provided we know four certain numbers that can be easily calculated for the sigma-words $C_k$, $D_k$, $C_{k+1}$ and $D_{k+1}$ (to be defined below), where $k=\lceil \log_2\ell \rceil$. Theorem~\ref{thm2} gives a recurrence relation for counting occurrences of patterns of the form $\tau_1\mn \tau_2$. In Section~\ref{section_6} we discuss occurrences of patterns of the form $\tau_1\mn \tau_2 \mn \cdots \mn \tau_k$, where the pattern $\tau_i$ does not overlap with the patterns $\tau_{i-1}$ and $\tau_{i+1}$ for $i=1,2,\ldots,k-1$. Finally, Section~\ref{9section7} deals with patterns of the form $[\tau_1\mn \tau_2\mn \cdots \mn \tau_k]$, $[\tau_1\mn \tau_2\mn \cdots \mn \tau_k)$ and $(\tau_1\mn \tau_2\mn \cdots \mn \tau_k]$ in Babson and Steingr\'{\i}msson notation, which means that we use "[x" in a pattern $p$ to indicate that in an occurrence of $p$, the letter corresponding to the $x$ must be the first letter of a word under consideration, whereas if we use ``$y]$'', we mean that the letter corresponding to $y$ must be the last (rightmost) letter in the word. 

\section{Preliminaries}\label{definitions}

In \cite{Evdok,Yab}, Evdokimov constructed chains of maximal length in the $n$-dimensional unit cube using the {\em sigma-sequence}. The sigma-sequence $w_{\sigma}$ was defined there by the following inductive scheme:

\centerline{$C_1=1$, \ \ \ \ $D_1=2$}

\centerline{$C_{k+1}=C_k1D_k$, \ \ \ \ $D_{k+1}=C_k2D_k$}

\centerline{$k=1, 2, \ldots$}

and $w_{\sigma}=\lim\limits_{ k \to \infty }{C_k}$. Thus, the initial letters of $w_{\sigma}$ are $11211221112212 \ldots$. We call the words $C_k$ the {\em sigma words}. The first four values of the sequence $\{C_k\}_{k\geq 1}$ are 1, 112, 1121122, 112112211122122.

In \cite{Kitaev} an equivalent definition of $w_{\sigma}$ was given: any natural number $n\neq 0$ can be presented  unambiguously as $n = 2^t(4s + \sigma)$, where $\sigma < 4$, and $t$ is the greatest natural number such that $2^t$ divides $n$. If $n$ runs through the natural numbers then $\sigma$ runs through some sequence consisting of 1 and 3. If we substitute 3 by 2 in this sequence, we get $w_{\sigma}$. 

In this paper we count occurrences of patterns in the sigma-words, which are particular initial subwords of $w_{\sigma}$. However, the challenging question is to find the number of occurrences of patterns or subwords in an arbitrary initial subword of $w_{\sigma}$, or more generally, in a subword of $w_{\sigma}$ starting in the position $i$ and ending in the position $j$.    

It turns out that for counting occurrences of certain patterns or subwords in $C_n$, one needs to know the number of occurrences of certain patterns in $D_n$. So, in the most cases, we give results for both $C_n$ and $D_n$. However, our main purpose is the words $C_n$ for $n\geq 1$, and in some propositions and examples we do not consider $D_n$.  

In what follows, we give initial values for the words $C_i$ and $D_i$:

\begin{tabular}{ll}
$C_1=1$ & $D_1=2$ \\[2mm]
$C_2=112$ & $D_2=122$ \\[2mm]
$C_3=1121122$ & $D_3=1122122$ \\[2mm]
$C_4=112112211122122$ & $D_4=112112221122122$ \\[2mm]
$C_5=1121122111221221112112221122122$ & $D_5=1121122111221222112112221122122$\\[2mm]
\end{tabular}

We now give some other definitions.

A {\em descent} (resp. {\em rise}) in a word $w=a_1a_2\ldots a_n$ is an $i$ such that $a_i > a_{i+1}$ (resp. $a_i < a_{i+1}$). It follows from the definitions that an occurrence of a descent (resp. rise) is an occurrence of the pattern $21$ (resp. $12$).

Let $c_n^{\tau}$ (resp. $d_n^{\tau}$) denote the number of occurrences of the pattern $\tau$ in $C_n$ (resp. $D_n$).

Suppose a word $W=AaB$, where $A$ and $B$ are some words of the same length, and $a$ is a letter. We define the {\em kernel of order} $k$ for the word $W$ to be the subword consisting of the $k-1$ rightmost letters of $A$, the letter $a$, and the $k-1$ leftmost letters of $B$. We denote it by $\KK_k(W)$. For example, $\KK_3(111211221)=12112$. If $|A|<k-1$ then we assume $\KK_k(W)=\epsilon$, that is, the kernel in this case is the empty word. Also, $m_k(\tau, W)$ denotes the number of occurrences of the pattern (or the word, or the subsequence depending on the context) $\tau$ in $\KK_k(W)$.   

We denote $x\mn x\mn \cdots \mn x$ ($\ell$ times) by $x^{\ell}$. Also, $\lceil a \rceil$ denotes the least natural number $b$ such that $a\leq b$.

\section{Patterns $1\mn 1\mn \cdots \mn 1$, $1\mn 2$ and $2\mn 1$}

It is easy to see that $|C_n|=|D_n|=2^{n}-1$. The following lemma gives the number of the letters 1 and 2 in $C_n$ and $D_n$.

\begin{lemma}\label{lemma1}
The number of 1s {\rm(}resp. 2s{\rm)} in $C_n$ is $2^{n-1}$ {\rm(}resp. $2^{n-1}-1${\rm)}. The number of 1s {\rm(}resp. 2s{\rm)} in $D_n$ is $2^{n-1}-1$ {\rm(}resp. $2^{n-1}${\rm)}. 
\end{lemma}

\begin{proof}
It is enough to find the number of 1s $c_n$ and $d_n$ in $C_n$ and $D_n$ respectively, since the number of 2s in $C_n$ and $D_n$ are obviously equal to $|C_n|-c_n$ and $|D_n|-d_n$ respectively.  

It is easy to see from the structure of $C_n$ and $D_n$ that 
\[
\left\{
\begin{array}{l}
c_n=c_{n-1}+d_{n-1}+1, \\[4pt]
d_n=c_{n-1}+d_{n-1},
\end{array}
\right.
\]
together with $c_1=1$ and $d_1=0$. The solution to this recurrence is $c_n=2^{n-1}$ and $d_n=2^{n-1}-1$.
\end{proof}

\begin{prop}\label{prop1}
The number occurrences of the subsequence $1^k$ {\rm(}resp. $2^k${\rm)} in $C_n$ is ${2^{n-1} \choose k}$ {\rm(}resp. ${2^{n-1}-1 \choose k}${\rm)}. Thus, the number of occurrences of the pattern $1^k$ in $C_n$ is equal to
$$c_n^{1^k}={2^{n-1} \choose k}+{2^{n-1}-1 \choose k}=\frac{2^n-k}{2^{n-1}-k}{2^{n-1}-1 \choose k}.$$ 
\end{prop}

\begin{proof}
From Lemma~\ref{lemma1}, there are $2^{n-1}$ (resp. $2^{n-1}-1$) occurrences of the letter $1$ (resp. 2) in $C_n$, and thus there are ${2^{n-1} \choose k}$ (resp. ${2^{n-1}-1 \choose k}$) occurrences of the subsequence $1^k$ (resp. $2^k$) there.   
\end{proof}

\begin{prop}\label{prop2}
We have that for all $n\geq 2$, $c_n^{1\mn 2}=d_n^{1\mn 2}=2\cdot4^{n-2}+(n-2)\cdot 2^{n-2}$, and $c_n^{2\mn 1}=d_n^{2\mn 1}=2\cdot 4^{n-2}-n\cdot 2^{n-2}$. 
\end{prop}

\begin{proof}
Let us first consider the pattern $1\mn2$. An occurrence of this pattern in $C_n=C_{n-1}1D_{n-1}$ is either inside $C_{n-1}$, or inside $D_{n-1}$, or the letter 1 is from the word $C_{n-1}1$, whereas the letter 2 is from the word $D_{n-1}$. Thus
$$c_n^{1\mn 2}=c_{n-1}^{1\mn2}+d_{n-1}^{1\mn2}+\mbox{\{ (the number of 1s in $C_{n-1}$) + 1\}} \cdot \mbox{\{ the number of 2s in $D_{n-1}$\}}.$$ 
Using the same considerations for $D_n=C_{n-1}2D_{n-1}$, one can get
$$d_n^{1\mn 2}=c_{n-1}^{1\mn2}+d_{n-1}^{1\mn2}+\mbox{\{ the number of 1s in $C_{n-1}$\}} \cdot \mbox{\{ (the number of 2s in $D_{n-1}$) + 1\}}.$$
The number of 1s and 2s in $C_{n-1}$ and $D_{n-1}$ is given in Lemma~\ref{lemma1}. So,
\[
\left\{
\begin{array}{l}
c_n^{1\mn2}=c_{n-1}^{1\mn2}+d_{n-1}^{1\mn2}+2^{n-2}\cdot(2^{n-2}+1)\\[6pt]
d_n^{1\mn2}=c_{n-1}^{1\mn2}+d_{n-1}^{1\mn2}+2^{n-2}\cdot(2^{n-2}+1)
\end{array}
\right.\Leftrightarrow
\]
\begin{equation}\label{1}
\left(
\begin{array}{l}
c_n^{1\mn2} \\[6pt] d_n^{1\mn2}
\end{array}
\right) =
\left(
\begin{array}{cc}
1 & 1  \\[6pt]
1 & 1  
\end{array}
\right)
\left(
\begin{array}{l}
c_{n-1}^{1\mn2} \\[6pt] d_{n-1}^{1\mn2}
\end{array}
\right)+
\left(
\begin{array}{l}
2^{n-2}\cdot(2^{n-2}+1) \\[6pt] 2^{n-2}\cdot(2^{n-2}+1)
\end{array}
\right)
\end{equation}
together with $c_2^{1\mn 2}=2$ and $d_2^{1\mn 2}=2$. Here, and several times in what follows, we need to solve recurrence relations of the form
$$x_n=Ax_{n-1}+b,$$
where $A$ is a matrix, and $x_n$, $x_{n-1}$ and $b$ are some vectors, where $b$ sometimes depends on $n$. We recall from linear algebra that such relations can be solved by diagonalization of the matrix $A$, that is, by writing $A=VDV^{-1}$, where $D$ is a diagonal matrix consisting of eigenvalues of $A$, and the columns of $V$ are eigenvectors of $A$. For example, if $A$ is a $2 \times 2$ matrix that consists of 1s, then we use 
$$
\left(
\begin{array}{cc}
1 & 1  \\
1 & 1  \\
\end{array}
\right)=
\left(
\begin{array}{rr}
 1 & 1 \\
-1 & 1 \\
\end{array}
\right)
\left(
\begin{array}{cc}
0 & 0  \\
0 & 2  \\
\end{array}
\right)
\left(
\begin{array}{rr}
1/2 & -1/2 \\
1/2 & 1/2 \\
\end{array}
\right)
$$ 
for computing powers of $A$, and thus for solving the recurrence relations. For the recurrence~\ref{1}, we get that for all $n\geq 2$, $c_n^{1\mn 2}=d_n^{1\mn 2}=2\cdot4^{n-2}+(n-2)\cdot 2^{n-2}$.

In the same manner, we can get that for the pattern $2\mn 1$,
\[
\left\{
\begin{array}{l}
c_n^{2\mn1}=c_{n-1}^{2\mn1}+d_{n-1}^{2\mn1}+2^{n-2}\cdot(2^{n-2}-1), \\[6pt]
d_n^{2\mn1}=c_{n-1}^{2\mn1}+d_{n-1}^{2\mn1}+2^{n-2}\cdot(2^{n-2}-1),
\end{array}
\right.
\]
together with $c_3^{2\mn 1}=2$ and $d_3^{2\mn 1}=2$. This gives, that for all $n\geq 2$, $c_n^{2\mn 1}=d_n^{2\mn 1}=2\cdot 4^{n-2}-n\cdot 2^{n-2}$.  
\end{proof}

Proposition~\ref{prop2} shows that asymptotically, the numbers of occurrences of the patterns, or the subsequences, $1\mn 2$ and $2\mn1$ in $C_n$ or $D_n$ are equal.

\section{Patterns without internal dashes}

Recall the definitions in Section~\ref{definitions}. 

\begin{thm}\label{thm1} Let $\tau=\tau_1\tau_2\ldots\tau_{\ell}$ be an arbitrary generalized pattern without internal dashes that consists of 1s and 2s. Suppose $k=\lceil \log_2\ell \rceil$, $a=m_{\ell}(\tau,D_k1C_k)$, and $b=m_{\ell}(\tau,D_k2C_k)$. Then for $n>k+1$, we have 
$$c_n^{\tau}=(a+b+c_{k+1}^{\tau}+d_{k+1}^{\tau})\cdot 2^{n-k-2}-b,$$
$$d_n^{\tau}=(a+b+c_{k+1}^{\tau}+d_{k+1}^{\tau})\cdot 2^{n-k-2}-a.$$
\end{thm}

\begin{proof}
Suppose $n>k+1$. In this case, $C_n=C_{n-1}1D_{n-1}=W_1\KK_{\ell}(D_k1C_k)W_2$, for some words $W_1$ and $W_2$ such that $|W_1|=|W_2|$. Because of the definition of the kernel $\KK_{\ell}(D_k1C_k)$, an occurrence of the pattern $\tau$ in $C_n$ is in either $C_{n-1}$, or $D_{n-1}$, or $\KK_k(D_k1C_k)$ (from the definitions $|C_{n-1} \cap \KK_k(D_k1C_k)| = |D_{n-1} \cap \KK_k(D_k1C_k)| = \ell -1$ and thus these intersections cannot be an occurrence of $\tau$). So,
$$c_n^{\tau}=c_{n-1}^{\tau}+d_{n-1}^{\tau} + a,$$ 
whereas in the same  way, we can obtain that 
$$d_n^{\tau}=c_{n-1}^{\tau}+d_{n-1}^{\tau} + b.$$
By solving these recurrence relations, we get the desirable. 
\end{proof}

In particular, Theorem~\ref{thm1} is valid for $\ell=1$, in which case the number of occurrences of $\tau$ in $C_n$ (or $D_n$) is the number of letters in $C_n$ (or $D_n$). Indeed, in this case, $k=0$, $a=b=c_1^1=d_1^1=1$, hence $c_n^1=d_n^1=2^n-1=|C_n|=|D_n|$. Also, as a corollary to Theorem~\ref{thm1} we have, that if $a=b=c_{k+1}^{\tau}=d_{k+1}^{\tau}=0$ for some pattern $\tau$, then this pattern never appears in sigma-sequence. 

All of the following examples are corollaries to Theorem~\ref{thm1}.

\begin{ex}\label{exam5} Suppose $\tau=12$. We have that $k=1$, $a=m_2(12,D_11C_1)=0$ and $b=m_2(12,D_12C_1)=0$. Besides, $c_2^{12}=1$ and $d_2^{12}=1$. Thus using Theorem~\ref{thm1}, for all $n>2$, $c_n^{12}=2^{n-2}$. So, the number of rises in $C_n$ is equal to $2^{n-2}$, for $n\geq 2$. 

If $\tau=21$, again $k=1$, but now $a=m_2(21,D_11C_1)=1$ and $b=m_2(21,D_12C_1)=1$. Besides, $c_3^{21}=1$ and $d_3^{21}=1$. From Theorem~\ref{thm1}, for all $n>3$, $c_n^{21}=2^{n-2}-1$, which shows that the number of descents in $C_n$ is one less than the number of rises. 

Since in both cases $a=b$, using the recurrences in Theorem~\ref{thm1}, we have that $c_n^{12}=d_n^{12}=2^{n-2}$, whereas $c_n^{21}=d_n^{21}=2^{n-2}-1$. \end{ex}

\begin{ex}\label{9example6} Suppose $\tau=112$. We have that $k=2$, $a=m_3(112,D_21C_2)=0$, and $b=m_3(112,D_22C_2)=0$. Besides, $c_3^{112}=2$ and $d_3^{112}=1$. Now, from Theorem~\ref{thm1}, we have that for all $n>3$, $c_n^{112}=d_n^{112}=3\cdot 2^{n-4}$.
\end{ex}

\begin{ex}\label{exam6} Suppose $\tau=221$. We have that $k=2$, $a=m_3(221,D_21C_2)=1$, and $b=m_3(221,D_22C_2)=1$. Besides, $c_3^{221}=0$ and $d_3^{221}=1$. Now, from Theorem~\ref{thm1}, we have that for all $n>3$, $c_n^{221}=d_n^{221}=3\cdot 2^{n-4}-1$. \end{ex}

\begin{ex} If $\tau=2212221$ then $k=3$, $a=m_7(221,D_31C_3)=0$, $b=m_7(221,D_32C_3)=1$, $c_4^{2212221}=0$, and $d_4^{2212221}=0$. Thus for $n\geq 4$, $c_n^{2212221}=2^{n-4}-1$. \end{ex}

\section{Patterns of the form $\tau_1\mn\tau_2$}\label{section_5}

\begin{thm}\label{thm2} Let $p=\tau_1\mn\tau_2$ be a generalized pattern such that $|\tau_1|=k_1$ and $|\tau_2|=k_2$. Suppose $k=\lceil \log_2(k_1+k_2-1)\rceil$. The following denote the number of occurrences of the subwords $\tau_1$ and $\tau_2$ in the certain kernels: $a_{\tau_1}=m_{k_1}(\tau_1,D_k1C_k)$, $a_{\tau_2}=m_{k_2}(\tau_2,D_k1C_k)$, $b_{\tau_1}=m_{k_1}(\tau_1,D_k2C_k)$, and $b_{\tau_2}=m_{k_2}(\tau_2,D_k2C_k)$. Also, let $r_1^a$ {\rm(}resp. $r_2^a$, $r_1^b$, $r_2^b${\rm)} denote the number of occurrences of overlapping subwords $\tau_1$ and $\tau_2$ in the word $D_k1C_k$ {\rm(}resp. $D_k1C_k$, $D_k2C_k$, $D_k2C_k${\rm)}, where $\tau_1\in \KK_{k_1}(D_k1C_k)$ and $\tau_2 \in C_k$ {\rm(}resp. $\tau_1 \in D_k$ and $\tau_2\in \KK_{k_2}(D_k1C_k)$, $\tau_1\in \KK_{k_1}(D_k2C_k)$ and $\tau_2 \in C_k$, $\tau_1 \in D_k$ and $\tau_2\in \KK_{k_2}(D_k2C_k)${\rm)}. Besides, we assume that we know $c_n^{\tau_i}$ and $d_n^{\tau_i}$ for $n>n_i$, $i=1,2$. Then for $n>\max(k+1,n_1+1,n_2+1)$, $c_n^{\tau}$ and $d_n^{\tau}$ are given by the following recurrence:
\[
\left(
\begin{array}{l}
c_n^{\tau} \\[6pt] d_n^{\tau}
\end{array}
\right) =
\left(
\begin{array}{cc}
1 & 1  \\[6pt]
1 & 1  
\end{array}
\right)
\left(
\begin{array}{l}
c_{n-1}^{\tau} \\[6pt] d_{n-1}^{\tau}
\end{array}
\right)+
\left(
\begin{array}{l}
\alpha_n \\[6pt] \beta_n
\end{array}
\right),
\] 
where
$$\alpha_n=(c_{n-1}^{\tau_1}+a_{\tau_1}-r_1^a)d^{\tau_2}_{n-1}+(a_{\tau_2}-r_2^a)c_{n-1}^{\tau_1}$$
and
$$\beta_n=(c_{n-1}^{\tau_1}+b_{\tau_1}-r_1^b)d^{\tau_2}_{n-1}+(b_{\tau_2}-r_2^b)c_{n-1}^{\tau_1}.$$
\end{thm}

\begin{proof} Suppose $n>\max(k+1,n_1+1,n_2+1)$. Let us find a recurrence for the number $c_n^{\tau}$ (one can use the same considerations for $d_n^{\tau}$).

An occurrence of the pattern $\tau$ in $C_n=C_{n-1}1D_{n-1}$ is either inside $C_{n-1}$, or inside $D_{n-1}$, or begins in $C_{n-1}$ or the letter 1 between $C_{n-1}$ and $D_{n-1}$ and ends in $D_{n-1}$ or the letter 1. The first two cases obviously give $c_{n-1}$ and $d_{n-1}$ occurrences of $\tau$. To count the contribution of the last to cases, we work with words instead of patterns. We do it to take in account the situations when $\tau_1$ or $\tau_2$ consists of copies of only one letter. In this case, we cannot count occurrence of these patterns separately, and then use this information, since, for instance, occurrences of the pattern $\tau_1=111$ are subwords 111 and 222 (the last one of these subwords we do not need), whereas occurrences of the pattern $\tau_1=222$ are not defined at all (222 is not a pattern). 

If an occurrence of $\tau_1\mn \tau_2$ does not entirely belong to $C_{n-1}$ or $D_{n-1}$ then we only have one of the following possibilities:
\begin{itemize}
\item[(a)] the subword $\tau_1$ entirely belongs to $C_{n-1}$ and the subword $\tau_2$ entirely belongs to $D_{n-1}$;
\item[(b)] the subword $\tau_1$ belongs entirely to $C_{n-1}$ and the subword $\tau_2$ belongs to the kernel $\KK_{k_2}(D_k1C_k)$, where $k=\lceil \log_2(k_1+k_2-1)\rceil$ is the least number that allow to control, in $C_n$ ($n>k$), overlapping occurrences of subwords $\tau_1$ and $\tau_2$ where $\tau_1$ is entirely from $C_{n-1}$ and $\tau_2 \in \KK_{k_2}(D_k1C_k)$;
\item[(c)] the subword $\tau_2$ belongs entirely to $D_{n-1}$ and the subword $\tau_1$ belongs to the kernel $\KK_{k_1}(D_k1C_k)$.
\end{itemize}

In (a) we obviously have $c_{n-1}^{\tau_1}\cdot d^{\tau_2}_{n-1}$ possibilities.

In (b) we have $c_{n-1}^{\tau_1}\cdot a_{\tau_2} - c_{n-1}^{\tau_1}\cdot r_2^a$ possibilities, since we need to subtract those occurrences of $\tau_1$ and $\tau_2$ that overlap.

Analogically to (b), in (c) we have  $d_{n-1}^{\tau_2}\cdot a_{\tau_1} - d_{n-1}^{\tau_2}\cdot r_1^a$ possibilities, which completes the proof.
\end{proof}

\begin{remark} For using Theorem~\ref{thm2}, one needs to know $c_n^{\tau}$ and $d_n^{\tau}$ for patterns $\tau$ without internal dashes. These numbers could be obtained by using Theorem~\ref{thm1}.
\end{remark}

The following corollary to Theorem~\ref{thm2} is straitforward to prove, using the fact that for non-overlapping patterns $\tau_1$ and $\tau_2$, $r_1^a=r_2^a=r_1^b=r_2^b=0$.

\begin{cor}\label{cor1} We make the same assumptions as those in Theorem~\ref{thm2}. Suppose additionally that the words $\tau_1$ and $\tau_2$ are not overlapping in the following sense: no one suffix of $\tau_1$ is a prefix of $\tau_2$. Then for $n>\max(k+1,n_1+1,n_2+1)$, $c_n^{\tau}$ and $d_n^{\tau}$ are given by the same recurrence as that in Theorem~\ref{thm2} with
$$\alpha_n=(c_{n-1}^{\tau_1}+a_{\tau_1})d^{\tau_2}_{n-1}+a_{\tau_2}c_{n-1}^{\tau_1}$$
and
$$\beta_n=(c_{n-1}^{\tau_1}+b_{\tau_1})d^{\tau_2}_{n-1}+b_{\tau_2}c_{n-1}^{\tau_1}.$$
\end{cor}

\begin{remark} Corollary~\ref{cor1} is valid under more weak assumptions, namely we only need the property of non-overlapping of the patterns $\tau_1$ and $\tau_2$ when one of them is in its kernel and the other one is not in its kernel. Example~\ref{hopelast} deals with the pattern $\tau$ that has overlapping blocks $\tau_1$ and $\tau_2$, but Corollary~\ref{cor1} can be applied. However, from practical point of view, checking the fact if two subwords are non-overlapping is more easy than considering the kernels and checking the non-overlapping of the subwords there.    
\end{remark}

\begin{ex} Suppose $\tau=12\mn 21$. We have that $|\tau_1|=|\tau_2|=2$. Now, in the statement of Theorem~\ref{thm2} we have that $k=2$, $a_{\tau_1}=0$, $a_{\tau_2}=1$, $b_{\tau_1}=0$ and $b_{\tau_2}=1$. Also, since there are no overlapping occurrences of the subwords $12$ and $21$ in $\KK_{3}(1221112)$ and $\KK_{3}(1222112)$, we have $r^a_1=0$, $r^a_2=0$, $r^b_1=0$ and $r^b_2=0$. Besides, from example~\ref{exam5}, $c_n^{12}=d_n^{12}=2^{n-2}$ and $c_n^{21}=d_n^{21}=2^{n-2}-1$. Thus, $\alpha_n=\beta_n=4^{n-3}$. Using the fact that $c_3^{12\mn21}=0$ and $d_3^{12\mn21}=1$, this allows us to get an explicit formula for $c_n^{12\mn 21}$ and $d_n^{12\mn 21}$ for $n>3$:
$$c_n^{12\mn 21}=d_n^{12\mn 21}=\frac{1}{2}4^{n-2}-3\cdot 2^{n-4}.$$
In particular $c_4^{12\mn21}=5$.     
\end{ex}

\begin{ex}\label{exam12} Suppose $\tau=1\mn 221$. We have that $|\tau_1|=1$ and $|\tau_2|=3$. Moreover, the words $\tau_1$ and $\tau_2$ are not overlapping, hence we can use Corollary~\ref{cor1}. We have that $k=2$, $a_{\tau_1}=1$, $a_{\tau_2}=1$, $b_{\tau_1}=0$ and $b_{\tau_2}=1$. From example~\ref{exam6}, $d_n^{221}=3\cdot 2^{n-4}-1$. Also, the number of occurrences of the letter 1 (the subword $\tau_1=1$) is given by Lemma~\ref{lemma1}: $c^1_n=2^{n-1}$. So, $\alpha_n=6\cdot 4^{n-4}+3\cdot 2^{n-5}-1$ and $\beta_n=6\cdot 4^{n-4}$. One can get now an explicit formula for $c_n^{1\mn 221}$ and $d_n^{1\mn 221}$ for $n>4$:
$$
\begin{array}{l}
c_n^{1\mn 221}=\frac{1}{2}4^{n-2}+27\cdot 2^{n-5}-n-7, \\[2mm]
d_n^{1\mn 221}=\frac{1}{2}4^{n-2}+21\cdot 2^{n-5}-8.
\end{array}
$$
In particular, $c_5^{1\mn 221}=47$.
\end{ex}

\begin{ex}\label{hopelast} Suppose $\tau=112\mn 21$. We have that $|\tau_1|=k_1=3$ and $|\tau_2|=k_2=2$. The other parameters in Theorem~\ref{thm2} are $k=3$, $a_{\tau_1}=0$, $a_{\tau_2}=1$, $b_{\tau_1}=0$, $b_{\tau_2}=1$, $r_1^a=r_2^a=r_1^b=r_2^b=0$. From Example~\ref{9example6}, for $n\geq 4$, $c_n^{112}=3\cdot 2^{n-4}$, and from Example~\ref{exam5}, $d_n^{21}=2^{n-2}-1$. Thus, in Theorem~\ref{thm2}, $\alpha_n=\beta_n=c_{n-1}^{112}(d_{n-1}^{21}+1)=3\cdot 4^{n-4}$. Now, we solve the recurrence relation from the theorem to get, that for $n>3$
$$c_n^{112\mn 21}=d_n^{112\mn 21}=\frac{3}{2}\cdot 4^{n-3}-2^{n-4}.$$
\end{ex}

\section{Counting occurrences of $\tau_1\mn \tau_2\mn \cdots \mn\tau_k$}\label{section_6}

In this section we study the number of occurrences of a pattern $\tau=\tau_1\mn \tau_2\mn \cdots \mn\tau_k$, where $\tau_i$ are patterns without internal dashes. We say that $\tau$ consists of $k$ blocks. We assume that for $i=1,2,\ldots,k-1$, the pattern $\tau_i$ does not overlap with the patterns $\tau_{i-1}$ and $\tau_{i+1}$. In this case we give a recurrence relation for the number of occurrences of $\tau$, provided we know the number of occurrences of certain patterns consisting of less than, or equal to, $k-1$ blocks, as well as $2k$ certain numbers which can be calculated by considering the words $D_{\ell}1C_{\ell}$ and $C_{\ell}2D_{\ell}$, where $\ell$ is the maximum number such that $\ell\leq \max_i\lceil \log_2|\tau_i| \rceil$. The cases of $k=1$ and $k=2$ are studied in the previous sections; they give the bases for our calculations. However, the case of overlapping patterns $\tau_i$ is not solved, and it remains as a challenging problem, since an answer to this problem gives the way to count occurrences of an arbitrary generalized pattern, or an arbitrary subsequence, in $\sigma$-words.   

\begin{thm}\label{mostgeneral}
Let $\tau=\tau_1\mn\tau_2\mn \cdots \mn\tau_k$ be a generalized pattern such that $|\tau_i|=k_i$ for $i=1,2,\ldots,k$. We assume that for $i=1,2,\ldots,k-1$, the subword $\tau_i$ does not overlap with the subwords $\tau_{i-1}$ and $\tau_{i+1}$ in the following sense: no one suffix of $\tau_{i-1}$ is a prefix of $\tau_i$ and no one suffix of $\tau_{i}$ is a prefix of $\tau_{i+1}$. Suppose $\ell_i=\lceil \log_2k_i\rceil$, $\ell=\max_i\ell_i$, and for the subwords $\tau_i$ we have $a_i=m_{k_i}(\tau_i,D_{\ell_i}1C_{\ell_i})$ and $b_i=m_{k_i}(\tau_i,D_{\ell_i}2C_{\ell_i})$, for $i=1,2,\ldots,k$. We assume that we know $c_{n-1}^{\tau_1\mn \cdots \mn \tau_i}$ and $d_{n-1}^{\tau_{i+1}\mn \cdots \mn \tau_k}$ for each $1\leq i \leq k-1$ and for all $n>n^{\star}$. Then for all $n>\max(\ell+1,n^{\star}+1)$, $c_n^{\tau}$ and $d_n^{\tau}$ are given by the following recurrence:
\[
\left(
\begin{array}{l}
c_n^{\tau} \\[6pt] d_n^{\tau}
\end{array}
\right) =
\left(
\begin{array}{cc}
1 & 1  \\[6pt]
1 & 1  
\end{array}
\right)
\left(
\begin{array}{l}
c_{n-1}^{\tau} \\[6pt] d_{n-1}^{\tau}
\end{array}
\right)+
\displaystyle\sum_{i=1}^{k-1}
\left(
\begin{array}{l}
c_{n-1}^{\tau_1\mn \cdots \mn \tau_i}\cdot d_{n-1}^{\tau_{i+1}\mn \cdots \mn \tau_k} \\[6pt]c_{n-1}^{\tau_1\mn \cdots \mn \tau_i}\cdot d_{n-1}^{\tau_{i+1}\mn \cdots \mn \tau_k}  
\end{array}
\right)+
\displaystyle\sum_{i=1}^{k}
\left(
\begin{array}{l}
a_i \cdot c_{n-1}^{\tau_1\mn \cdots \mn \tau_{i-1}}\cdot d_{n-1}^{\tau_{i+1}\mn \cdots \mn \tau_k} \\[6pt] b_i \cdot c_{n-1}^{\tau_1\mn \cdots \mn \tau_{i-1}}\cdot d_{n-1}^{\tau_{i+1}\mn \cdots \mn \tau_k}   
\end{array}
\right).
\] 
\end{thm}

\begin{proof}
We consider only $c_n^{\tau}$, since the same arguments can be applied to $d_n^{\tau}$. We use the considerations similar to those in Theorem~\ref{thm2}.

An occurrence of the pattern $\tau$ in $C_n=C_{n-1}1D_{n-1}$ can be entirely in $C_n$ or $D_n$. The first term counts such occurrences. Otherwise, we have two possibilities: either the letter 1 between the words $C_{n-1}$ and $D_{n-1}$ does not belong to an occurrence of $\tau$, or it does do it, in which case there exist $i$ (exactly one) such that the subword $\tau_i$ occurs in its kernel. The first sum in the statement is obviously responsible for the first of this cases, whereas the second sum is responsible for the second case (in the last case we use the fact that subwords $\tau_i$ are not overlapping).   
\end{proof}

As a corollary to Theorem~\ref{mostgeneral}, we have Corollary~\ref{cor1}.

The following example is another corollary to Theorem~\ref{mostgeneral}.

\begin{ex} Suppose $\tau=2-1-221$, that is, $\tau_1=2$, $\tau_2=1$ and $\tau_3=221$. So, parameters in Theorem~\ref{mostgeneral} are the following: $k_1=k_2=1$, $k_3=3$, $\ell_1=\ell_2=1$, $\ell_3=2$, $\ell=2$. From $D_11C_1=211$ we obtain $a_1=0$, $a_2=1$. From $D_21C_2=1221112$ we obtain $a_3=1$. From $D_12C_1=221$ we get $b_1=1$, $b_2=0$. From $D_22C_2=1222112$ we get $b_3=1$. Besides, from Proposition~\ref{prop2}, Examples~\ref{exam6} and~\ref{exam12}, we have
\[
\begin{array}{l}
c_n^{\tau_1\mn \tau_2}=c_n^{2-1}=2\cdot4^{n-2}-n\cdot 2^{n-2}, \mbox{for $n>1$};\\[3mm]
d_n^{\tau_3}=d_n^{221}=3\cdot 2^{n-4}-1, \mbox{for $n>3$}; \\[3mm]
d_n^{\tau_2\mn \tau_3}=d_n^{1\mn 221}=\frac{1}{2}\cdot 4^{n-2}+21\cdot 2^{n-5}-8, \mbox{for $n>4$}.
\end{array}
\]
Also, the number of occurrences of the subword $\tau_1=2$ in $C_n$ is given by Proposition~\ref{prop1}: $c_n^{\tau_1}=c_n^2=2^{n-1}-1$. So, the number of occurrences of the pattern $\tau$ in $C_n$ and $D_n$, for $n>5$, satisfies the following recurrence relation:
\[
\left(
\begin{array}{l}
c_n^{\tau} \\[6pt] d_n^{\tau}
\end{array}
\right) =
\left(
\begin{array}{cc}
1 & 1  \\[6pt]
1 & 1  
\end{array}
\right)
\left(
\begin{array}{l}
c_{n-1}^{\tau} \\[6pt] d_{n-1}^{\tau}
\end{array}
\right)+
\left(
\begin{array}{l}
\frac{5}{1024}8^n+\frac{25-3n}{256}4^n-\frac{171}{64}2^n+9 \\[6pt]
\frac{5}{1024}8^n+\frac{21-3n}{256}4^n-2^{n+1}
\end{array}
\right),
\] 
with initial conditions $c_5^{\tau}=70$ and $d_5^{\tau}=74$.
\end{ex}

\section{Patterns of the form $[\tau_1\mn \tau_2\mn \cdots \mn \tau_k]$, $[\tau_1\mn \tau_2\mn \cdots \mn \tau_k)$ and $(\tau_1\mn \tau_2\mn \cdots \mn \tau_k]$}\label{9section7}

We recall that according to Babson and Steingr\'{\i}msson notation for generalized patterns, if we use "[" in a pattern, for example if we write $p=[1\mn 2)$, we indicate that in an occurrence of $p$, the letter corresponding to the 1 must be the first letter of a word under consideration, whereas if we write, say, $p=(1\mn 2]$, then the letter corresponding to 2 must be the last (rightmost) letter of the word. 

In the theorems of this section, we assume that we can find the numbers $c_n^{\tau_1\mn \tau_2\mn \cdots \mn \tau_k}$ and $d_n^{\tau_1\mn \tau_2\mn \cdots \mn \tau_k}$ for any patterns $\tau_i$, $i=1,2,\ldots,k$, without internal dashes. For certain special cases, these numbers can be obtained using the theorems of Sections~\ref{section_5} and~\ref{section_6}. 

\begin{thm}\label{9theorem15} Suppose $\tau_1$ and $\tau_2$ are two patterns without internal dashes such that $|\tau_1|=k_1$ and $|\tau_2|=k_2$. Also, suppose $\ell_1=\log_2(k_1+1)$, $\ell_2=\log_2(k_2+1)$ and $\ell=\log_2(k_1+k_2+1)$. 
Let $a(\tau_1,\tau_2)$ be the number of overlapping subwords $\tau_1$ and $\tau_2$ in $C_{\ell}$ such that 
$\tau_1$ occurs as $k_1$ leftmost letters of $C_{\ell}$; $b(\tau_1,\tau_2)$ is the number of overlapping subwords $\tau_1$ and $\tau_2$ in $C_{\ell}$ such that $\tau_2$ occurs as $k_2$ rightmost letters of $C_{\ell}$. We assume that we know $c_n^{\tau_i}$ and $d_n^{\tau_i}$ for $i=1,2$ and for all $n>n^{\star}$.
\begin{itemize}
\item[i.] For $n\geq \max(\ell_1,n^{\star})$,
\[
c_n^{[\tau_1\mn \tau_2)}=
\left\{
\begin{array}{ll}
0, & \mbox{if $C_{\ell_1}$ does not begin with $\tau_1$},\\[6pt]
c_n^{\tau_2}-a(\tau_1,\tau_2), & \mbox{otherwise}.
\end{array}
\right.
\]
\item[ii.] For $n\geq \max(\ell_2,n^{\star})$,
\[
c_n^{(\tau_1\mn \tau_2]}=
\left\{
\begin{array}{ll}
0, & \mbox{if $C_{\ell_2}$ does not end with $\tau_2$},\\[6pt]
c_n^{\tau_1}-b(\tau_1,\tau_2), & \mbox{otherwise}.
\end{array}
\right.
\]
\item[iii.] For $n\geq \ell$,
\[
c_n^{[\tau_1\mn \tau_2]}=
\left\{
\begin{array}{ll}
0, & \mbox{if $C_{\ell}$ does not begin with $\tau_1$ or end with $\tau_2$},\\[6pt]
1, & \mbox{otherwise}.
\end{array}
\right.
\]
\item[iv.] For $n\geq \max(\ell_1,n^{\star})$,
\[
d_n^{[\tau_1\mn \tau_2)}=
\left\{
\begin{array}{ll}
0, & \mbox{if $D_{\ell_1}$ does not begin with $\tau_1$},\\[6pt]
d_n^{\tau_2}-a(\tau_1,\tau_2), & \mbox{otherwise}.
\end{array}
\right.
\]
\item[v.] For $n\geq \max(\ell_2,n^{\star})$,
\[
d_n^{(\tau_1\mn \tau_2]}=
\left\{
\begin{array}{ll}
0, & \mbox{if $D_{\ell_2}$ does not end with $\tau_2$},\\[6pt]
d_n^{\tau_1}-b(\tau_1,\tau_2), & \mbox{otherwise}.
\end{array}
\right.
\]
\item[vi.] For $n\geq \ell$,
\[
d_n^{[\tau_1\mn \tau_2]}=
\left\{
\begin{array}{ll}
0, & \mbox{if $D_{\ell}$ does not begin with $\tau_1$ or end with $\tau_2$},\\[6pt]
1, & \mbox{otherwise}.
\end{array}
\right.
\]
\end{itemize}
\end{thm}

\begin{proof} We prove case i, all the other cases are then easy to see.

Clearly, if $C_{\ell_1}$ does not begin with $\tau_1$ then $C_n$ does not begin with $\tau_1$ for all $n\geq \ell_1$, which means that $c_n^{[\tau_1\mn \tau_2)}=0$ in this case. Otherwise, to count occurrences of the pattern $[\tau_1\mn \tau_2)$ is the same as to find the number of occurrences of the pattern $\tau_2$ in $C_n$ and then subtract the number of such occurrences of $\tau_2$ that begin from the $i$-th letter of $C_n$, where $1\leq i\leq k_1$. 
\end{proof}

The following two examples are corollaries to Theorem~\ref{9theorem15}.

\begin{ex} Suppose we have the patterns $\sigma_1=[1122-21211)$ and $\sigma_2=(21221-12]$. From Theorem~\ref{9theorem15}, $c_n^{\sigma_1}=d_n^{\sigma_1}=0$, since $C_3$ does not begin with $1122$ ($\ell_1=3$). Also, $c_n^{\sigma_2}=d_n^{\sigma_2}=0$, since $C_3$ does not end with $12$ ($\ell_2=3$).
\end{ex}

\begin{ex} Suppose $\tau=[112\mn 21)$. We have that $k_1=3$, $\ell_1=2$ and $C_2$ begins with the subword $112$. Besides, $a(112,21)=1$ and, from Example~\ref{exam5}, $c_n^{21}=d_n^{21}=2^{n-2}-1$. Theorem~\ref{9theorem15} now gives, that for $n> 3$, we have $c_n^{[112\mn 21)}= c_n^{\tau_2}-a(\tau_1,\tau_2)= 2^{n-2} - 2$. \end{ex}

The following theorem is straitforward to prove using the assumptions concerning non-overlapping of certain subwords.

\begin{thm}\label{9thm18} Let $\{\tau_1, \tau_2, \ldots ,\tau_k\}$ be a set of generalized patterns without internal dashes. Suppose $|\tau_1|=s_1$, $|\tau_k|=s_k$, $\ell_1=\log_2(s_1+1)$ and $\ell_k=\log_2(s_k+1)$. Also, $\ell=\max(\ell_1,\ell_k)$.

\begin{itemize}
\item[i.] With the assumption that the subword $\tau_1$ does not overlap with the subword $\tau_2$, that is, no one suffix of $\tau_1$ is a prefix of $\tau_2$, we have
\begin{itemize}
\item[(a)] 
\[
c_n^{[\tau_1\mn \tau_2\mn \cdots \mn \tau_k)}=
\left\{
\begin{array}{ll}
0, & \mbox{if $C_{\ell_1}$ does not begin with $\tau_1$},\\[6pt]
c_n^{\tau_2\mn \tau_3\mn \cdots \mn \tau_k}, & \mbox{otherwise}.
\end{array}
\right.
\]
\item[(b)]
\[
d_n^{[\tau_1\mn \tau_2\mn \cdots \mn \tau_k)}=
\left\{
\begin{array}{ll}
0, & \mbox{if $D_{\ell_1}$ does not begin with $\tau_1$},\\[6pt]
d_n^{\tau_2\mn\tau_3\mn \cdots \mn \tau_k}, & \mbox{otherwise}.
\end{array}
\right.
\]
\end{itemize}
\item[ii.] With assumption that the subword $\tau_{k-1}$ does not overlap with the subword $\tau_k$, that is, no one suffix of $\tau_{k-1}$ is a prefix of $\tau_k$, we have 
\begin{itemize}
\item[(a)] 
\[
c_n^{(\tau_1\mn \tau_2\mn \cdots \mn \tau_k]}=
\left\{
\begin{array}{ll}
0, & \mbox{if $C_{\ell_k}$ does not end with $\tau_k$},\\[6pt]
c_n^{\tau_1\mn\tau_2\mn \cdots \mn \tau_{k-1}}, & \mbox{otherwise}.
\end{array}
\right.
\]
\item[(b)]
\[
d_n^{(\tau_1\mn \tau_2\mn \cdots \mn \tau_k]}=
\left\{
\begin{array}{ll}
0, & \mbox{if $D_{\ell_k}$ does not end with $\tau_k$},\\[6pt]
d_n^{\tau_1\mn\tau_2\mn \cdots \mn \tau_{k-1}}, & \mbox{otherwise}.
\end{array}
\right.
\]
\end{itemize}

\item[iii.] With the assumption that the subword $\tau_1$ does not overlap with the subword $\tau_2$, and the subword $\tau_{k-1}$ does not overlap with the subword $\tau_k$, we have 
\begin{itemize}
\item[(a)] 
\[
c_n^{[\tau_1\mn \tau_2\mn \cdots \mn \tau_k]}=
\left\{
\begin{array}{ll}
0, & \mbox{if $C_{\ell}$ does not begin with $\tau_1$ or does not end with $\tau_k$},\\[6pt]
c_n^{\tau_2\mn\tau_3\mn \cdots \mn \tau_{k-1}}, & \mbox{otherwise}.
\end{array}
\right.
\]
\item[(b)]
\[
d_n^{[\tau_1\mn \tau_2\mn \cdots \mn \tau_k]}=
\left\{
\begin{array}{ll}
0, & \mbox{if $D_{\ell}$ does not begin with $\tau_1$ or does not end with $\tau_k$},\\[6pt]
d_n^{\tau_2\mn\tau_3\mn \cdots \mn \tau_{k-1}}, & \mbox{otherwise}.
\end{array}
\right.
\]
\end{itemize}
\end{itemize}
\end{thm}

The following example is a corollary to Theorem~\ref{9thm18}.

\begin{ex} Suppose $\tau=[112\mn 1\mn 221\mn 22]$. The parameters of Theorem~\ref{9thm18} are $k_1=3$, $k_2=2$, $\ell_1=2$, $\ell_2=2$, $\ell=2$. $C_3$ begins with the subword $112$ and ends with the subword $22$. Thus by Theorem~\ref{9thm18} and Example~\ref{exam12}, $c_n^{[112\mn 1\mn 221\mn 22]}=c_n^{1\mn 221}=\frac{1}{2}4^{n-2}+27\cdot 2^{n-5}-n-7$.
\end{ex}


\begin{thebibliography}{10}

\bibitem[BabStein]{BabStein} Babson E., Steingr\'{\i}msson E.: Generalized permutation patterns and a classification of the Mahonian statistics, S\'em. Lothar. Combin. {\bf 44} (2000), Art. B44b, 18 pp.
\bibitem[Burstein]{Burstein} Burstein A., Enumeration of words with forbidden patterns, Ph.D.
thesis, University of Pennsylvania, (1998).
\bibitem[BurMans1]{BurMans1} Burstein A., Mansour T.: Words restricted by patterns with at
most 2 distinct letters, \emph{Electronic J. of Combinatorics}, to appear (2002).
\bibitem[BurMans2]{BurMans2} Burstein A., Mansour T.: Words restricted by $3$-letter
generalized multipermutation patterns, preprint CO/0112281.
\bibitem[BurMans3]{BurMans3} Burstein A., Mansour T.: Counting occurrences of some subword
patterns, preprint CO/0204320.
\bibitem[Claes]{Claes} A. Claesson: Generalised Pattern Avoidance, European J. Combin. {\bf 22} (2001), no. 7, 961--971.
\bibitem[Evdok]{Evdok} Evdokimov A. A.: On the Maximal Chain Length of an Unit $n$-dimensional Cube, Maths Notes {\bf 6}, No. 3 (1969), 309--319. (Russian)
\bibitem[Evdokimov]{Evdokimov} Private communication (2001).
\bibitem[GelbOlm]{GelbOlm} Gelbaum B., Olmsted J.: {\em Counterexamples in Analysis}, Holden-day, San Francisco, London, Amsterdam, (1964).
\bibitem[Kitaev]{Kitaev} Kitaev S., There are no iterated morphisms that define the Arshon sequence and the sigma-sequence, to appear J. Automata, Languages and Combinatorics (2002).
\bibitem[KitMans1]{KitMans1} Kitaev S., Mansour T.: Counting the occurrences of generalized patterns in words generated by a morphism, preprint CO/0210170.
\bibitem[KitMans2]{KitMans2} Kitaev S., Mansour T.: The Peano curve and counting occurrences of some patterns, preprint CO/0210268.
\bibitem[Knuth]{Knuth} Knuth D. E.: {\em The Art of Computer Programming}, 2nd ed. Addison Wesley, Reading, MA, (1973).
\bibitem[Lothaire]{Lothaire} Lothaire M.: {\em Combinatorics on Words}, Encyclopedia of Mathematics, Vol. {\bf 17}, Addison-Wesley (1986). Reprinted in the {\em Cambridge Mathematical Library}, Cambridge University Press, Cambridge UK (1997).
\bibitem[Salomaa]{Salomaa} Salomaa A.: {\em Jewels of Formal Language Theory}, Computer Science Press (1981).
\bibitem[SimSch]{SimSch} Simion R., Schmidt F.: Restricted permutations, European J. Combin. {\bf 6}, no. 4 (1985), 383--406.
\bibitem[Yab]{Yab} Yablonsky S. V.: {\em Discrete mathematics and mathematical problems of cybernetics}, Nauka, Vol.~{\bf 1}, Moscow (1974), 112--116. (Russian)
\end{thebibliography}
\end{document}